\documentclass[a4paper,leqno]{article}
\usepackage{amsmath}
\usepackage{latexsym}
\usepackage{amssymb}
\usepackage{amsthm}
\newcommand{\dis}{\displaystyle}
\newcommand{\m}{\mathbb}
\newtheorem{thm}{Theorem}
\newtheorem{df}[thm]{Definition}
\newtheorem{prop}[thm]{Proposition}
\newtheorem{lem}[thm]{Lemma}						
\newtheorem{cor}[thm]{Corollary}
\newtheorem{rem}[thm]{Remark}
\newtheorem{ex}[thm]{Example}

\newcommand{\ze}[1]{\zeta _{MT,{#1}} }

\newcommand{\unb}{\underbrace}
\newcommand{\nm}[2]{ \{ #1 \}^{#2} }
\newcommand{\bb}[2]{ #1 _1,#1 _2,\dots ,#1 _{#2} }
\newcommand{\bbp}[2]{ #1 _1 + \cdots + #1 _{#2} }
\newcommand{\bbg}[2]{0< #1 _1 < #1 _2 < \cdots < #1 _{#2} }
\newcommand{\bbpara}[3]{ {#1}^{#2 _1}_1 {#1}^{#2 _2}_2 \cdots {#1}^{{#2}_{#3}}_{#3} }
\DeclareMathOperator{\Li}{Li}
\DeclareMathOperator{\La}{\mathcal{L}}
\DeclareMathOperator{\ii}{\int ^\infty _0 \int ^\infty _0 \cdots \int ^\infty _0}

\title{ On analogues of the Arakawa-Kaneko zeta functions of Mordell-Tornheim type}
\author{\Large Takuma Ito }
\date{}

\begin{document}

\maketitle

\begin{abstract}
In this paper, we construct certain analogues of the Arakawa-Kaneko zeta functions. We prove functional relations between these functions and the Mordell-Tornheim multiple zeta functions. Furthermore we give some formulas among Mordell-Tornheim multiple zeta values as their applications.
\end{abstract}

\section{Introduction}
Let $\m{Z}$ be the rational integer ring, $\m{N}$ the set of natural numbers, $\m{N}_0=\m{N} \cup \{ 0 \}$, $\m{Q}$ the rational number field and $\m{C}$ the complex number field. We denote $n$ repetitions of $m$ by $\nm{m}{n}$ for $m,n \in \m{N}$. 

Arakawa and Kaneko $\cite{A1}$ introduced the ``Arakawa-Kaneko zeta function'' defined by
\begin{align}\label{zz}
\xi (\bb{k}{r}; s) = \frac{1}{\Gamma (s)} \int ^\infty _0 \frac{t^{s-1} }{e^t -1} \Li _{\bb{k}{r}}(1-e^{-t} ) dt 
\end{align}
for $(\bb{k}{r}) \in \m{N}^r$ and $s \in \m{C}$ with $\Re(s)>0$, where $\Li _{\bb{k}{r}} (z)$ is the polylogarithm defined by 
\begin{align*}
\Li _{\bb{k}{r}} (z) =\sum _{\bbg{m}{r}} \frac{z^{m_r}}{\bbpara{m}{k}{r}} \quad (z \in \m{C},\ \lvert z \rvert <1).
\end{align*}
\\
\\
When $r=1$, $\xi(k;s)$ is also denoted by $\xi_k(s)$. They proved that for $m \in \m{N}_0$, $\xi (\nm{1}{r-1},k;m+1)$ can be written in terms of multiple zeta values (MZVs) in \cite[Theorem 9]{A1}. 

On the other hand, Matsumoto defined the ``Mordell-Tornheim $r$-ple zeta function'' by 
\begin{align}\label{mt}
&\quad \ze{r} (\bb{s}{r};s_{r+1}) \\
&=\sum ^\infty_{\bb{m}{r}= 1} \frac{1}{\bbpara{m}{s}{r} \left(\sum ^r _{j=1} m_j \right)^{s_{r+1}}} \quad (s_i \in \m{C}, \Re(rs_i+s_{r+1} ) > r), \nonumber
\end{align}
and proved that this function can be continued meromorphically to the whole $\m{C}^{r+1}$-space in \cite{M1} and \cite{M2}. This zeta function in the double sum case was first studied by Tornheim \cite{T} for the values at positive integers in 1950s. He gave some evaluation formulas for $\ze{2} (k_1,k_2;k_3)$ for $k_1, k_2, k_3 \in \m{N}$. Mordell \cite{Mor} independently proved that $\ze{2} (k,k;k)  \pi ^{-3k} \in \m{Q}$ for all even $k \geq 2$. Tsumura \cite[Theorem 4.5]{Tsumura1} and Nakamura \cite[Theorem 1]{Naka1} showed certain functional relations among the Mordell-Tornheim double zeta functions and the Riemann zeta functions. 

In this paper, for $\Bbbk \in \m{N}^r$, we first define the function 
\[
\xi _{MT} (\Bbbk ;s)=\frac{1}{\Gamma(s)}\int ^\infty _0 \frac{t^{s-1}}{e^t -1} \prod^r _{j=1} \Li _{k_j}(1-e^{-t}) dt \quad (s \in \m{C},\ \Re(s)>0)
\]
which can be regarded as an analogue of the Arakawa-Kaneko zeta function of    
Mordell-Tornheim type (see Definition $\ref{XMT}$).  
We construct functional relations between $\xi_{MT}(\Bbbk;s)$  and the Mordell-Tornheim multiple zeta functions (see Theorem $\ref{MR1}$).  For example, 
\begin{align*}
&\quad \zeta(2) ^2 \zeta(s) -2\zeta(2) \xi _{MT}(2; s) + \xi _{MT} (2,2;s) \\
&= \zeta  _{MT,3} (2,2,0;s)  +2s \zeta _{MT,3} (2,1,0;s+1) +s(s+1) \zeta _{MT,3} (1,1,0;s+2). \nonumber
\end{align*}
This can be proved by the method similar to the proof of \cite[Theorem 8]{A1}. 

Secondly, we show certain relation formulas among Mordell-Tornheim multiple zeta values (see Corollary $\ref{MR2}$). For example, 
\begin{align*}
&\quad \zeta(2) ^2 \zeta(m+1) -2\zeta(2) \frac{1}{m!}\ze{m+1}(2,\nm{1}{m};1) 
+ \frac{1}{m!}\ze{m+2}(2,2,\nm{1}{m};1) \\
&= \ze{3} (2,2,0;m+1) +2(m+1) \ze{3} (2,1,0;m+2) \nonumber \\
&\quad +(m+1)(m+2) \ze{3} (1,1,0;m+3) \quad (m \in \m{N},\ m \geq 3),  \nonumber 
\end{align*}
and
\[
\ze{3}(2,1,1;1)= 2\zeta(2)\zeta(3)-\zeta(5). 
\]
Lastly, we consider a generalization of main results (see Theorem $\ref{MR4}$). 
\section{Preliminaries} 
We first construct a Mordell-Tornheim type analogue of $\xi(\bb{k}{r};s)$ and continue it analytically to an entire function. We define $\{ C^\Bbbk _{m,MT} \}$ by 
\begin{align}\label{aqs}
\frac{\prod ^r _{j=1} \Li _{k_j}(1-e^{-t})}{e^t-1} =\sum ^\infty _{m=0} C^{\Bbbk}_{m,MT} \frac{t^m}{m!}
\end{align}
for $\Bbbk =(\bb{k}{r}) \in \m{Z} ^r$. These are generalizations of poly-Bernoulli numbers $\{ C ^{(k)} _m \}$ defined by 
\[
\frac{\Li _k (1-e^{-t})}{e^t -1} = \sum ^\infty _{m=0} C^{(k)} _m \frac{t^m}{m!}
\]
for $k \in \m{Z}$ (see \cite{A1}). Since $\Li _k(1-e^{-t})=O(t) \ (t \to 0)$ and $\Li_k(1-e^{-t})=O(t) \ (t \to \infty)$ for all $k \in \m{N}$, we can define the following function. 
\begin{df}\label{XMT}
For $\Bbbk=(\bb{k}{r}) \in \m{N}^r$ and $s \in \m{C}$ with $ \Re(s) >1-r$, let
\begin{align}\label{ximt}
\xi _{MT} (\Bbbk ;s)=\frac{1}{\Gamma(s)}\int ^\infty _0 \frac{t^{s-1}}{e^t -1} \prod^r _{j=1} \Li _{k_j}(1-e^{-t}) dt,
\end{align}
where $\Gamma(s)$ is the gamma function. 
\end{df}
The integral on the right-hand side of $(\ref{ximt})$ converges absolutely uniformly in the region $\Re(s) > 1-r$. When $r=1$, $\xi _{MT} (\Bbbk ;s)=\xi_k (s)$ holds for $\Bbbk=(k) \in \m{N}$. 

\begin{thm}\label{AC1}
For $\Bbbk  =(\bb{k}{r}) \in \m{N}^r$, the function $\xi _{MT} (\Bbbk ;s)$ can be continued analytically to an entire function, and satisfies
\begin{align}
\xi _{MT} (\Bbbk ;-m) = (-1)^m C^\Bbbk _{m,MT} \quad (m \in \m{N}_0). 
\end{align}
\begin{proof}
Let
\begin{align}
A(\Bbbk ;s) &= \int _C \frac{t^{s-1}}{e^t -1} \prod ^r _{j=1} \Li _{k_j} (1-e^{-t} ) dt  \nonumber \\
&= \left( e^{2\pi \sqrt{-1} s } -1 \right) \int ^\infty _\varepsilon \frac{t^{s-1} }{e^t -1} \prod ^r _{j=1} \Li _{k_j} (1-e^{-t} )dt \nonumber \\ 
&\quad + \int _{C_\varepsilon} \frac{t^{s-1} }{e^t -1} \prod ^r _{j=1} \Li _{k_j} (1-e^{-t} ) dt \quad (s \in \m{C}), \nonumber
\end{align}
where $C$ is the contour which is the path consisting of the real axis (top side), a circle $C_\varepsilon$ around the origin of radius $\varepsilon$ (sufficiently small), and the positive real axis (bottom side). 
Since the integrand has no singularity on $C$ and the contour integral converges absolutely for all $s \in \m{C}$, we can see that $A(\Bbbk ;s)$ is entire. 
Suppose $\Re(s) > 1-r$, then the second integral tends to $0$ as $\varepsilon \to 0$. Therefore we have
\begin{align*}
\xi _{MT} (\Bbbk ;s) = \frac{1}{\left( e^{2 \pi \sqrt{-1} s}-1 \right)\Gamma(s) }A(\Bbbk ;s).
\end{align*}
Since $\xi _{MT} (\Bbbk ;s)$ is holomorphic for $\Re(s) >1-r$, this function has no singularity at any positive integer. Therefore this gives the analytic continuation of $\xi _{MT} (\Bbbk ;s)$ to an entire function. Let $s= -m$ for $m\in \m{N}_0$. Using $(\ref{aqs})$, we have
\begin{align*}
\xi _{MT} (\Bbbk ;-m) &= \frac{(-1)^m m!}{2\pi \sqrt{-1} } A(\Bbbk ;-m) \nonumber \\
&= \frac{(-1)^m m!}{2 \pi \sqrt{-1} } \int _{C_\varepsilon} t^{-m-1} \sum ^\infty _{n=0} 
C^\Bbbk _{n,MT} \frac{t^n}{n!} dt \nonumber \\
&=  (-1)^m C^\Bbbk _{m,MT}. 
\end{align*}
This completes the proof.
\end{proof}
\end{thm}
Secondly, we show a relation between the Mordell-Tornheim multiple zeta values and $\xi_{MT}(\Bbbk ;m+1)$ for $m \in \m{N}_0$. For this aim, we consider the following function and give a lemma. 
\begin{df} For $\Bbbk =(\bb{k}{r+1}) \in \m{N} ^{r} \times \m{N}_0$ and $z \in \m{C}$ with $\lvert z \rvert <1$, let
\begin{align}\label{mtpoly}
\La_\Bbbk (z) =\sum ^\infty _{\bb{m}{r} =1} \frac{z^{\sum ^r _{j=1} m_j}}{\bbpara{m}{k}{r} \left(\sum ^r _{j=1} m_j \right)^{k_{r+1}} }. 
\end{align}
\end{df}
Under the above condition, the sum on the right-hand side of $\eqref{mtpoly}$ converges absolutely uniformly. We note that $\La_\Bbbk (z) = \Li_{k_1+k_2}(z)$ holds for $r=1$ and $\Bbbk=(k_1,k_2)$. By direct calculation, we have
\begin{lem}\label{AAAA} For $\Bbbk =(\bb{k}{r},k_{r+1} ) \in \m{N}^{r+1}$ and $z \in \m{C}$ with $\lvert z \rvert <1$, 
\[
\frac{d}{dz} \La _\Bbbk (z)=
\begin{cases}
\dis \frac{1}{z}  \La _{\underline{\Bbbk}^{(r+1)} } (z) & (k_{r+1} \geq 2), \\
\\
\dis \frac{1}{z}  \prod ^r_{j=1} \Li _{k_j} (z) & (k_{r+1} =1),
\end{cases}
\]
where $\underline{\Bbbk}^{(r+1)}=(\bb{k}{r},k_{r+1} -1)$.
\end{lem}
Using Lemma $\ref{AAAA}$ and calculating directly, we obtain
\begin{prop}\label{spvalofxi} For $\Bbbk_1=(\bb{k}{r},0) \in \m{N}^r \times \m{N}_0,\ \Bbbk=(\bb{k}{r})$ and $m \m \in \m{N}_0$, 
\begin{align}
\xi_{MT}(\Bbbk ;m+1) = \frac{1}{m!} \ze{m+r}(\bb{k}{r},\nm{1}{m};1).  \nonumber
\end{align}
\end{prop}
We can recover \cite[Corollary 4.2 and Theorem 4.4]{Hof} as follows. 
\begin{cor}\label{MTval}
For $m \in \m{N}_0$, 
\[
\ze{m+1}(\nm{1}{m+1};1)=(m+1)!\zeta(m+2).
\]
\begin{proof}
By $\xi_1(s)=s\zeta(s+1)$ and Proposition $\ref{spvalofxi}$, we obtain the assertion.
\end{proof}
\end{cor}
\section{Main results}
In this section, we give main results. We first prepare the following lemma which is necessary 
to show the first and second main results.
\begin{lem}\label{int2} For $s_j \in \m{C}$ with $ \Re(s_j) >0\ (2 \leq j \leq r)$ and $\Re(s_{r+1}) >r$, 
\begin {align}
&\quad \ze{r} (0,s_2,\dots ,s_{r};s_{r+1} ) = \\
&\frac{1}{\prod ^{r+1}_{j =2} \Gamma(s_j)}  \underbrace{ \ii }_{r}  
 \frac{\prod^{r+1} _{j=2} t^{s_j -1} _j }{(e^{t_{r+1}}-1)\prod ^{r}_{j =2} (e^{t_j+t_{r+1}}-1)} dt_2\cdots dt_{r} dt_{r+1}. \nonumber
\end{align}
\begin{proof}
Using the well-known relation
\begin{align*}
m ^{-s}  = \frac{1}{\Gamma(s)} \int ^\infty _0 t^{s-1} e^{-m t} dt \quad (m \in \m{N},\ s \in \m{C},\ \Re(s)>0 )
\end{align*}
for $s_j \in \m{C}$ with $\Re(s_j) >0\ (2 \leq j \leq r)$ and $\Re(s_{r+1})>r$, we have
\begin{align*}
&\quad \prod ^{r+1} _{j=2} \Gamma(s_j) \times \ze{r}(0,s_2,\dots ,s_r;s_{r+1}) \nonumber \\
&= \sum ^\infty _{\bb{m}{r}=1} \prod ^r _{j=2} \left( \int ^\infty _0 t^{s_j-1} _j e^{-m_j t} dt_j \right) \nonumber \\
&\quad \times  \left( \int ^\infty _0 t^{s_{r+1}-1}_{r+1} e^{-\left(\sum ^r _{j=1} m_j \right)t_{r+1} } dt_{r+1} \right) \nonumber \\
&= \sum ^\infty _{\bb{m}{r}=1} \unb{\ii}_{r} dt_2 dt_3 \cdots dt_{r+1} \nonumber \\ 
&\quad \times  \left( \prod ^{r+1} _{j=2} t^{s_j -1}_j \right) \left( e^{-m_1 t_{r+1}} \right) \left( \prod ^r _{j=2}  e^{-m_j(t_j + t_{r+1} )} \right) \nonumber \\
&= \unb{\ii}_{r}  \frac{\prod ^{r+1}_{j=2} t^{s_j-1}_j}{\left(e^{t_{r+1}-1}-1 \right)\prod ^{r}_{j=2} \left( e^{t_j+t_{r+1}}-1 \right)}  dt_2 dt_3 \cdots dt_{r+1}.  \nonumber
\end{align*}
Changing the order of summation and integration is justified by absolutely convergence.  Therefore we complete the proof.
\end{proof}
\end{lem}
Using Lemma $\ref{int2}$, we have the first main result as follows. 
\begin{thm}\label{MR1}
For $r \in \m{N}$ and $s \in \m{C}$,  
\begin{align*}
& \quad \sum ^{r-1} _{j=0} \binom{r-1}{j} (-1)^j \zeta(2)^{r-1-j} \xi _{MT}  (\nm{2}{j} ;s)  \nonumber \\
&=\sum ^{r-1} _{j=0} \binom{r-1}{j} (s)_j \zeta  _{MT,r} (\nm{2}{r-1-j}, \nm{1}{j} ,0;s+j).
\end{align*}
\end{thm}
\begin{proof}
We first assume $r \geq2$. For $s \in \m{C}$ with $\Re(s) >0$, let
\begin{align*}
J_{MT,r} (s)=\underbrace{\int ^\infty _0 \cdots \int ^\infty _0}_{r} dt_1 dt_2 \cdots dt_r \frac{t^{s-1} _r}{e^{t_r}-1} \prod ^{r-1} _{j=1} \frac{t_j+t_r}{e^{t_j+t_r}-1}. 
\end{align*}
Using 
\[
\frac{\partial}{\partial t_j} \Li _2(1-e^{-t_j-t_r})= \frac{t_j+t_r}{e^{t_j+t_r}-1} \quad (1 \leq j \leq r-1), 
\]
we have
\begin{align*}
 J_{MT,r} (s) 
&= \sum ^{r-1} _{j=0} \binom{r-1}{j} (-1)^{r-1-j} \zeta (2) ^{r-1-j} \Gamma(s) \xi _{MT} (\nm{2}{j} ;s).
\end{align*}
On the other hand, by Lemma $\ref{int2}$ and 
\begin{align*}
\ze{r} (\dots ,\stackrel{i}{\Check{s_i}} ,\dots ,\stackrel{j}{\Check{s_j}},\dots ;s_{r+1} ) = \ze{r} (\dots ,\stackrel{i}{\Check{s_j}},\dots ,\stackrel{j}{\Check{s_i}},\dots ;s_{r+1} )
\end{align*}
for $1\leq i \leq j \leq r$, we have
\begin{align*}
J_{MT,r} (s)
&= \sum ^{r-1} _{j=0} \binom{r-1}{j} \Gamma(s+j) \ze{r} (\nm{2}{r-1-j},\nm{1}{j},0;s+j).
\end{align*}
By the analytic continuation, we obtain the desired identity in the case $r \geq 2$. 
When $r=1$, it holds trivially. Therefore we complete the proof.
\end{proof}
By Theorem $\ref{MR1}$ and Proposition $\ref{spvalofxi}$, we immediately obtain the second main result as follows. 
\begin{cor}\label{MR2}
For $r \in \m{N}$ and $m \in \m{N}_0$, 
\begin{align*}
&\quad \sum ^{r-1} _{j=0} \binom{r-1}{j} \frac{(-1)^j \zeta(2)^{r-1-j} }{m!} \zeta _{MT,j+m}  ( \nm{2}{j},\nm{1}{m} ;1)  \\
&=\sum ^{r-1} _{j=0} \binom{r-1}{j} (m+1)_j  \zeta  _{MT,r} (\nm{2}{r-1-j},\nm{1}{j} ,0;m+1+j). \nonumber
\end{align*}
\end{cor} 
Next, in order to evaluate $\ze{2k+1}(2,\nm{1}{2k} ;1)$, we quote \cite[(75)]{BBB}:
\begin{align}\label{dzv}
\zeta(a,b) &= \frac{1}{2} \left\{ \left((-1)^b \binom{M}{a}-1 \right) \zeta(M) +(1+(-1)^b )\zeta(a) \zeta(b) \right\} \\
&\quad +(-1)^{b+1} \sum ^{ (M-3)\slash 2} _{k =1} \left\{ \binom{2k}{a-1} + \binom{2k}{b-1} \right\}\zeta(2k+1)  \zeta(M-2k-1), \nonumber 
\end{align}
{\it where $a,b \in \m{N}$ with $a,b \geq2$ and $M=a+b \equiv 1 \pmod 2$.}
\begin{rem}
{\rm We note that $\eqref{dzv}$ also holds for $a=1$ providing we remove the term containing $\zeta(1)$.}  
\end{rem}
Combining $\eqref{dzv}$ and Corollary $\ref{MR2}$ in the case $r=2$, we have the third main result as follows.
\begin{prop}\label{MR3}
For $k \in \m{N}$, 
\begin{align*}
\ze{2k+1}(2,\nm{1}{2k};1)
&= (2k)!\biggl\{ \zeta(2)\zeta(2k+1) -\frac{1}{2}(2k^2+k-2)\zeta(2k+3) \\
&\quad +\sum ^{k-1}_{n=1} (2k+1-2n)\zeta(2n+1)\zeta(2k+2-2n)  \biggl\} .
\end{align*}
\end{prop}
\begin{ex}
\begin{align*}
\ze{3}(2,1,1;1)&= 2\zeta(2)\zeta(3)-\zeta(5), \\
\ze{5}(2,1,1,1,1;1)&= 4!\left\{ \zeta(2)\zeta(5)+3\zeta(3)\zeta(4) -4\zeta(7) \right\}.
\end{align*}
\end{ex}
These results correspond to \cite[Theorems 6, 8, 9 and Corollary 11]{A1}. Results in \cite{A1} are relations between $\xi (\bb{k}{r} ;s)$ and multiple zeta functions or MZVs. On the other hand, our results are relations between $\xi _{MT}(\Bbbk ;s)$ and Mordell-Tornheim multiple zeta functions or Mordell-Tornheim multiple zeta values.  
\section{A generalization of the function $\xi _{MT} (\Bbbk ;s)$}

In this section, we consider a certain generalization of the function $\xi _{MT} (\Bbbk ;s)$ and aim to generalize Theorem $\ref{MR1}$. 

By the definition $\eqref{mtpoly}$, for $\Bbbk =(\bb{k}{r},k_{r+1}) \in \m{N} ^r \times \m{N}_0$, we have
\begin{gather}\label{LemLa}
\La _\Bbbk (1-e^{-t}) = 
\begin{cases}
O(t^l) & \text{if $k_{r+1}=0$ and $l=\sharp \{j \ \vert \ k_j =1 \} \geq 1$} ,\\
O(1) & \text{otherwise} \quad (t \to \infty)
\end{cases}
\intertext{and} 
\label{ordofLa}
\La _\Bbbk (1-e^{-t}) = O(t^r)\ \left(t \rightarrow 0 \right).
\end{gather}
Using $\eqref{LemLa}$ and $\eqref{ordofLa}$, we can define the following function. 
\begin{df}\label{genexi}
For $\bb{r}{g} \in \m{N},\ \Bbbk _{i} =(\bb{k^{(i)}}{r_i},k ^{(i)}_{{r_i}+1} ) \in \m{N}^{r_i} \times \m{N}_0$ and $s \in \m{C}$ with $\Re(s) > 1-\sum ^g _{i=1} r_i$, let
\begin{align}\label{defofgenexi}
\xi _{MT,g} (\bb{\Bbbk}{g} ;s) = \frac{1}{\Gamma(s)} \int ^\infty _0 \frac{t^{s-1}}{e^t -1} \prod ^g _{i=1}\La _{\Bbbk _i} (1-e^{-t}) dt.
\end{align}
\end{df}
The integral on the right-hand side of $\eqref{defofgenexi}$ converges absolutely uniformly in the region $\Re(s)> 1-\sum ^g _{i=1} r_i$. Further we note that 
\[
\xi_{MT,1}(\Bbbk_1 ;s)=\xi_{MT}(\Bbbk ;s)  
\]
for $\Bbbk _1=(\bb{k}{r},0),\ \Bbbk =(\bb{k}{r})$ and $s \in \m{C}$. Therefore we can see that Definition $\ref{genexi}$ is a generalization of the function $\xi _{MT}(\Bbbk ;s)$. By the same method as in the proof of Theorem $\ref{AC1}$, we have
\begin{thm}\label{AC2}
For $g,\bb{r}{g} \in \m{N}$ and $\Bbbk_i \in \m{N}^{r_i} \times \m{N}_0 \ (1 \leq i \leq g)$, the function $\xi _{MT,g} (\bb{\Bbbk}{g} ;s)$ can be continued analytically to an entire function. 
\end{thm}
By the same method as in the proof of Theorem $\ref{MR1}$, we obtain
\begin{thm}\label{MR4} For $N \in \m{N},\ \mathbf{r} = (\bb{r}{N-1}) \in \m{N}^{N-1}$ and $ s \in \m{C}$, 
\begin{align}\label{gene1}
&\quad \sum ^{N-1} _{n=0}  (-1)^{n} \sum _{\substack{J \subset I_{N-1} \\ \sharp J =n}} \left\{ \prod _{j \in I_{N-1} \setminus J} \ze{r_j} (\nm{1}{r_j} ;1) \right\} \xi _{MT,n} (\{ \mathbf{1}_{r_j+1}\ \vert \ j \in  J \} ;s) \\
&= \sum ^{\mathrm{wt}(\mathbf{r} )} _{n=0} (s)_n \sum _{\substack{\bbp{i}{N-1}=n \\ r_l \geq i_l \geq 0}} \left\{ \prod ^{N-1} _{j=1} \binom{r_j}{i_j} (r_j-i_j) ! \right\}  \nonumber \\
&\quad \times  \ze{N} \left(r_1-i_1+1,r_2-i_2+1,\dots ,r_{N-1}-i_{N-1}+1,0 ; s+n \right), \nonumber
\end{align}
where $\mathbf{1}_{r_j+1} = (\nm{1}{r_j+1} ) \in \m{N}^{r_j+1},\ \xi _{MT,N} ( \emptyset ;s)=\zeta(s), \ I_{N-1}=\{ 1,2,\dots ,N-1 \}$ and $\mathrm{wt}(\mathbf{r})=\sum ^{N-1}_{i=1}r_i$.
\begin{proof}
We first assume $N \geq 2$ and define the function $J_{MT, \mathbf{r}}(s)$ by
\begin{align*}
J_{MT,\mathbf{r}} (s) &= \underbrace{\ii}_{N} dt_1 dt_2 \cdots  dt_N \\
&\quad \times  \frac{t_{N} ^{s-1} }{e^{t_{N} } -1} \prod ^{N-1} _{j=1} \frac{(t_j+t_{N} )^{r_j} }{e^{t_j+t_{N} }-1} \quad (s \in \m{C},\ \Re(s)>0)
\end{align*}
for $\mathbf{r}=(\bb{r}{N-1}) \in \m{N}^{N-1}$. It follows from Lemma $\ref{AAAA}$ that 
\[
\frac{\partial}{\partial t_j} \La _{\mathbf{1} _{r_j +1} } (1-e^{-t_j-t_N }) = \frac{(t_j+t_N)^{r_j} }{e^{t_j+t_N} -1}. 
\]
Therefore we have
\begin{align*}
J_{MT,\mathbf{r}}(s)
&= \sum ^{N-1} _{n=0}  (-1)^{n} \sum _{\substack{J \subset I_{N-1} \\ \sharp J =n}} \left\{ \prod _{j \in I_{N-1} \setminus J} \ze{r_j} (\nm{1}{r_j} ;1) \right\} \\
&\quad \times \Gamma(s) \xi _{MT,n} (\{ \mathbf{1}_{r_j+1}\ \vert \ j \in  J \} ;s)
\end{align*}
for $\Re(s)>1$. On the other hand, by Lemma $\ref{int2}$, we have
\begin{align*}
J_{MT,\mathbf{r}}(s) 
&= \sum ^{\mathrm{wt}(\mathbf{r} )} _{n=0} \Gamma\left( s+n \right) \sum _{\substack{\bbp{i}{N-1}=n \\ r_l \geq i_l \geq 0}} \left\{ \prod ^{N-1} _{j=1} \binom{r_j}{i_j} (r_j -i_j)! \right\}  \\
&\quad \times \ze{N} \left(r_1-i_1+1,r_2-i_2+1,\dots ,r_{N-1}-i_{N-1}+1,0 ; s+n \right)
\end{align*}
for $\Re(s)>N$. By the analytic continuation, we obtain $\eqref{gene1}$ for all $s \in \m{C}$ when $N \geq 2$. When $N=1$, $\eqref{gene1}$ holds obviously. Therefore the proof is completed. 
\end{proof}
\end{thm}
\begin{rem}
{\rm 
In particular, 
Theorem $\ref{MR4}$ in the case $(N, \mathbf{r})=(r, \mathbf{1}_{r-1})$ coincides with Theorem $\ref{MR1}$. Hence we can see that Theorem $\ref{MR4}$ is a generalization of Theorem $\ref{MR1}$.  
}
\end{rem}
We have not obtained the values of $\xi _{MT,g} (\bb{\Bbbk}{g} ;m+1)$ for $m \in \m{N}_0$. But we have a certain lemma as follows. 
\begin{lem}\label{lll} For $g, \bb{r}{g} \in \m{N}$,  
\begin{align*}
&\quad \sum ^g _{j=1} r_j ! \xi _{MT,g-1} (\mathbf{1}_{r_1+1},\dots,\mathbf{1}_{r_{j-1}+1},\mathbf{1}_{r_{j+1}+1},\dots ,\mathbf{1}_{r_g+1} ;r_j+1) \\
&=\prod ^g _{j=1}  \ze{r_j}(\nm{1}{r_j} ;1).
\end{align*}
\end{lem}
\begin{rem}
{\rm In particular, combining Corollary $\ref{MTval}$, Theorem $\ref{MR4}$ in the case $N=2$ and Lemma $\ref{lll}$ in the case $g=2$, we have the Euler decomposition (cf. \cite{A1}). }
\begin{align*}
\zeta (k+1) \zeta(r+1) &= \sum ^k _{m=0} \binom{r+m}{r}\zeta(k+1-m, r+1+m)  \\
&\quad + \sum ^r_{n=0} \binom{k+n}{k} \zeta(r+1-n, k+1+n) \quad (r,k \in \m{N}). \nonumber 
\end{align*}
\end{rem}

\section*{acknowledgments}
The author thanks Professor Hirofumi Tsumura for useful advice and pointing out some mistakes and unsuitable expressions.

Takuma Ito \\
Department of Mathematics and Information Science, Tokyo Metropolitan \\ University, 1-1, Minami-Ohsawa, Hachioji, Tokyo 192-0397 Japan \\
e-mail: sugakunotakuma.ito@gmail.com
\end{document}